\documentclass{elsart}
\usepackage{graphics}
\usepackage{graphicx}
\usepackage{epsfig}

\usepackage{amssymb}
\usepackage{enumerate}
\usepackage[latin1]{inputenc}
\usepackage[OT1]{fontenc}
\usepackage{color}

\usepackage{amsmath}
\usepackage{amsfonts}
\usepackage{makeidx}

\newtheorem{remark}[thm]{Remark}

 \newtheorem{proof}[thm]{Proof}

\newtheorem{algorithm12}[thm]{Algorithm}

\newtheorem{corrolar}[thm]{Corrolar}

\newcommand{\pa}{\partial}
\newcommand{\be}{\begin{equation}}
\newcommand{\ee}{\end{equation}}
\newcommand{\ba}{\begin{eqnarray}}
\newcommand{\ea}{\end{eqnarray}}
\newcommand{\nn}{\nonumber}
\newcommand{\la}{\label} 
\newcommand{\ep}{\Delta t}
\newcommand{\T}{{\mathcal T}} 
\newcommand{\dt}{{\Delta t}} 
\newcommand{\bp}{{\bf p}}
\newcommand{\bq}{{\bf q}}

\newcommand{\bv}{{\bf v}}

\begin{document}

\begin{frontmatter}



\title{Multiscale methods for Levitron Problems: Theory and Applications.}


\author[label1]{J\"urgen Geiser\corauthref{cor1}}
\ead{juergen.geiser@uni-greifswald.de}
\ead[url]{http://www.mathematik.hu-berlin.de/$\sim$geiser/}
\corauth[cor1]{J\"urgen Geiser}
\address[label1]{EMA University of Greifswald, Department of Physics, Felix-Hausdorff-Str. 6, D-17489 Greifswald, Germany}

\begin{abstract}
In this paper, we describe a multiscale model
based on magneto-static traps of neutral atoms or ion traps.
The idea is to levitate a magnetic spinning top 
in the air repelled by a base magnet.
 
For such a problem, we have to deal with different
time and spatial scales and we propose a novel
 splitting method for solving the levitron problem, see \cite{geiser_2012}.

We focus on the multiscale problem, which we obtain 
by coupling the kinetic $T$ and the potential $U$ part
of our equation.
The kinetic and potential parts, can be seen as generators of flows, see \cite{dullin2004}.

The main problem is based on the accurate computation
of the Hamiltonian equation and we propose a novel higher order
splitting scheme to obtain stable states near the relative equilibrium.
To improve the splitting scheme we apply a novel method so 
called MPE (multiproduct expansion method), see \cite{chin_2011}, which
include higher order extrapolation schemes.

In numerical studies, we discuss the stability near this 
relative equilibrium with our improved time-integrators. 
Best results are obtained by iterative and extrapolated Verlet 
schemes in comparison to higher order explicit Runge-Kutta schemes.
Experiments are applied to a magnetic top in an axisymmetric magnetic field 
(i.e. the Levitron) and we discuss the future applications
to quantum computations.

\end{abstract}


\end{frontmatter}

% main text

\section{Introduction}
\label{intro}

We are motivated to simulate a Levitron, which is a magnetic spinning top and
can levitate in a magnetic field. 
The main problem of such a nonlinear problem is to achieve a stability
for the calculation of  the critical splint rate.
While the stability of Levitrons are discussed in the work of \cite{dull98}
and their dynamics in \cite{gans97}, we concentrate on improving the
standard time-integrator schemes for the Hamiltonian systems.
It is important to derive stable numerical schemes with high accuracy
to compute the non-dissipative equation of motions, which are at least
higher order symplectic integrators.
Here, we apply geometric integrators based on the Str\"omer-Verlet method
with extrapolation methods, see \cite{chin_2011}.
While we have symplectic schemes, we preserve the underlying
physics of our Levitron, i.e. reversibility, symplecticity, volume preservation
and conservation of the first integrals, see \cite{hair2003}.

For the numerical studies, we propose novel splitting 
schemes and analyze their behavior.
We deal with a standard Verlet integrator and improve its accuracy with
iterative and extrapolation ideas. Such a Hamiltonian splitting method, 
can be seen as geometric integrator and saves computational time
while decoupling the full equation system.

The paper is organized as follows.
A mathematical model based on a multiscale problem of the Levitron 
is introduced in Section \ref{model}.
The splitting method is used as a solver method to decouple the
multiscale equations to more simpler equations given in the 
kinetics and the potential of the models is described in 
Section \ref{splitting}.
The improvement of the splitting strategies based on extrapolation 
schemes are discussed in Section \ref{extrapolation}.
The numerical experiments and their description of our used methods
are described in Section \ref{numerical}.
Finally the conclusions and on overview for our next works are discussed
in Section \ref{conclusion}.

\section{Mathematical model}
\label{model}

The Levitron is described on the base of rigid body theory. With the convention of Goldstein \cite{goldstein81} for the Euler angles the angular velocity $\omega_{\phi}$ is along the $z$-axis of the system, $\omega_{\theta}$ along the line of nodes and $\omega_{\psi}$ along the $z'$-axis.
Finally the kinetic energy can be written as 

\be
T=\frac{1}{2}\left[m(\dot x^2+\dot y^2+\dot z^2) + A( \dot \theta ^2+\dot \phi ^2 \sin ^2 \theta) + C(\dot \psi +\dot \phi \cos \theta)^2 \right]
\ee
The potential energy $U$ is given by the sum of the gravitational energy and the interaction potential of the Levitron in the magnetic field of the base plate:

\be 
U=mgz-\mu (\sin \psi \sin \theta \frac{\Phi}{x} + \cos \psi \sin \theta \frac{\Phi}{y} + \cos \theta \frac{\Phi}{z})
\ee  

with $\mu$ as the magnetic moment of the top and $\Phi$ the magneto-static potential. Following Gans \cite{gans97} we uses the potential of a ring dipole as approximation for a magnetized plane with a centered unmagnetized hole. Furthermore we introduced a nondimensionalization for the variables and the magneto-static potential:

\be
\Psi=\frac{Z}{(1+Z^2)^{3/2}}-(X^2+Y^2)\frac{3}{4}\frac{(2Z^2-3)Z}{(1+Z^2)^{7/2}}
\ee

Lengths were scaled by the radius R of the base plane, mass were measured in units of $m$ and energy in units of $mgh$. Therefore the one time unit is $\sqrt{R/g}$.

Knowing the kinetic and the potential energy, we can formulate the the Lagrangian as:
\begin{align}
 \mathcal{L}= & T(\dot x^2,\dot y^2,\dot z^2) - U(x,y,z) \notag \\
= &\frac{1}{2} \left[ m(\dot x^2 + \dot y^2 + \dot z^2)+ A(\dot\theta^2+\dot \phi^2 \sin^2\theta) + C(\dot\psi+\dot\phi\cos\theta)^2  \right] \notag \\
 & + \mu \left[\sin\theta \sin\psi B_x+ \sin\theta \cos\psi B_y + \cos\theta B_z \right] - mgz 
\end{align}
Furthermore, the Hamiltonian $\mathcal H$ can be calculated as 
\begin{align}
 \mathcal{H}	=&\vec{\dot q}^T \vec{p} - \mathcal{L} \notag \\
	=&\frac{1}{2m}\left(p_x^2+p_y^2+p_z^2\right) + \frac{p_\theta ^2}{2A} + \frac{p_\psi ^2}{2C} + \frac{\left(p_\phi -p_\psi \cos\theta \right)^2}					{2A\sin ^2\theta} \notag \\
	& +\mu \left[\sin\theta \sin\psi \frac{\partial \varphi}{\partial x}+ \sin\theta \cos\psi \frac{\partial \varphi}{\partial y} + \cos\theta 						\frac{\partial 		\varphi}{\partial z} \right] + mgz \, .
\end{align}
where 
\begin{align}
 \vec{\dot q}= &(\dot x;\dot y;\dot z;\dot \theta;\dot \psi;\dot \phi) \notag\\
 				=&\left(\frac{p_x}{m};\frac{p_y}{m},\frac{p_z}{m};\frac{p_\theta}{A};\frac{p_\psi}{C}-\frac{p_\phi \cos \theta - p_\psi \cos ^2 \theta}{A\sin ^2\theta}; \frac{p_\phi - p_\psi \cos \theta}{A\sin ^2\theta} \right)  
\end{align}
and ${\bf p}$ is given as:
\begin{align}
 p_x&=\frac{\partial \mathcal{L}}{\partial \dot x}=m \dot x \\
 p_y&=\frac{\partial \mathcal{L}}{\partial \dot y}=m \dot y \\
 p_z&=\frac{\partial \mathcal{L}}{\partial \dot z}=m \dot z \\
 p_\theta &=\frac{\partial \mathcal{L}}{\partial \dot \theta}=A\dot\theta \\
 p_\psi &=\frac{\partial \mathcal{L}}{\partial \dot \psi}=C(\dot\psi + \dot \phi \cos \theta) \\
 p_\phi &=\frac{\partial \mathcal{L}}{\partial \dot \phi}= A\dot\phi \sin^2 \theta + C(\dot \psi + \dot\phi \cos\theta)\cos\theta 
\end{align}

The Hamiltonian $\mathcal H$ is calculated as 
\begin{align}
 \mathcal{H}	=&\vec{\dot q}^T \vec{p} - \mathcal{L} \notag \\
	=&\frac{1}{2m}\left(p_x^2+p_y^2+p_z^2\right) + \frac{p_\theta ^2}{2A} + \frac{p_\psi ^2}{2C} + \frac{\left(p_\phi -p_\psi \cos\theta \right)^2}					{2A\sin ^2\theta} \notag \\
	& +\mu \left[\sin\theta \sin\psi \frac{\partial \varphi}{\partial x}+ \sin\theta \cos\psi \frac{\partial \varphi}{\partial y} + \cos\theta 						\frac{\partial 		\varphi}{\partial z} \right] + mgz \, .
\end{align}

In the next section, we discuss the time-integrator methods to solve
our differential equations.

\section{Splitting Methods}
\label{splitting}

The evolution of the dynamical variable $u(\bq,\bp)$ (including $\bq$ and $\bp$ themselves)
is given by the
Poisson bracket,
\begin{eqnarray}
\label{splitt_1}
\pa_tu(\bq,\bp)=
                 \Bigl(
		          {\frac{\partial u}{\partial \bq}}\cdot
                  {\frac{\partial H}{\partial \bp}}
				 -{\frac{\partial u}{\partial \bp}}\cdot
                  {\frac{\partial H}{\partial \bq}}
				                    \Bigr)=(A + B)u(\bq,\bp).
\end{eqnarray}

$A$ and $B$ are Lie operators, or vector fields
\be
A= \frac{\partial H}{\partial \bp} \cdot\frac{\pa}{\pa\bq} \qquad B= -  \frac{\partial H}{\partial \bq} \cdot\frac{\pa}{\pa\bp}
\la{shop} 
\ee

The transfer to the operators are given in the following description.

The exponential operators $\e^{\ep A}$ and $\e^{\ep B}$ are then just shift operators,
with $\T_2(\ep)$ is a symmetric second order splitting method:
\be
\T_{2, VV}(\ep)=\e^{(\ep/2) B}\e^{\ep A} \e^{(\ep/2) B}.
\la{str2}
\ee
and corresponds to the velocity form of the Verlet algorithm (VV).

\begin{figure}[htbp]
\centering	
	\begin{minipage}[b]{0.47\textwidth}
	\centering 
	\includegraphics[width=\textwidth]{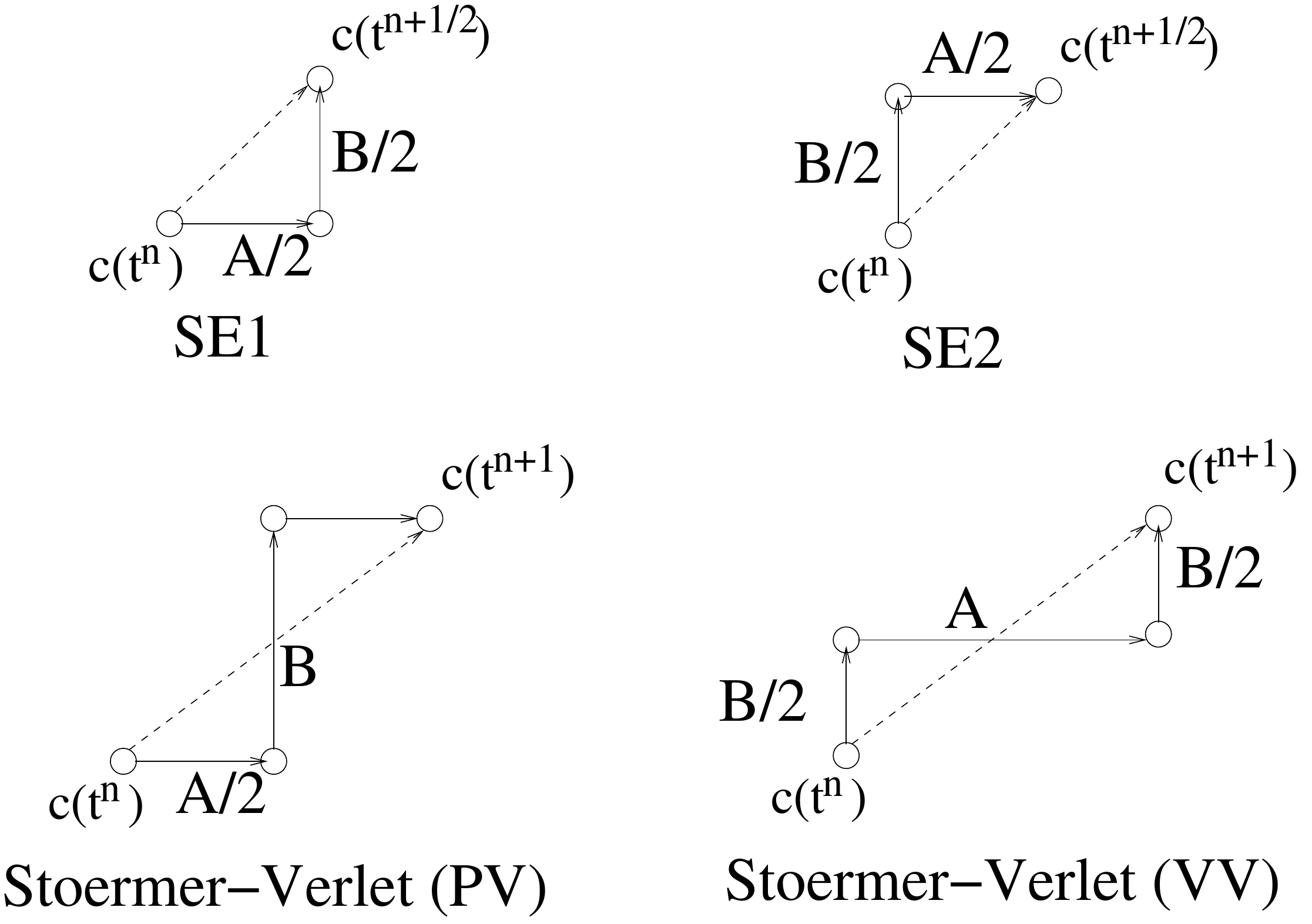}
	\end{minipage}
\caption{Symplectic Splitting Methods.} 
 \label{verlet}
\end{figure}

Further the splitting scheme:
\be
\T_{2, PV}(\ep)=\e^{(\ep/2) A}\e^{\ep B} \e^{(\ep/2) A}.
\la{str2}
\ee
and corresponds to the position-form of the Verlet algorithm (PV). \\
In the literature, see \cite{hair2003}, they are also known as symplectic 
splitting methods, see Figure \ref{verlet}.

The symplectic St\"ormer-Verlet or leap-frog algorithm in the
notation $\T_{2, VV}(\ep) = SE2(\ep/2) \circ SE1(\ep/2)$ is given in the 
following algorithm \ref{algo_2}.

\begin{algorithm12}
\label{algo_2}

We start with $(\bq_0, \bp_0)^t = (\bq(t^{n}), \bp(t^{n}))^t $:
\begin{eqnarray}
(\bq_1, \bp_1)^t = \e^{\ep/2 B} (\bq_0, \bp_0)^t & = & ( I - \frac{1}{2} \ep \sum_i  \frac{\partial H}{\partial \bq}(\bp_i,  \bq_i)  \frac{\partial}{\partial \bp_i}) (\bq_0, \bp_0)^t ,
\end{eqnarray}

\begin{eqnarray}
(\bq_2, \bv_2)^t = \e^{\ep A} (\bq_1, \bv_1)^t & = &  ( I + \ep \sum_i  \frac{\partial H}{\partial \bp}(\bp_i,  \bq_i)  \frac{\partial}{\partial \bq_i}) (\bq_1, \bp_1)^t , 
\end{eqnarray}

\begin{eqnarray}
(\bq_3, \bv_3)^t =  \e^{\ep/2 B} (\bq_2, \bp_2)^t & = & ( I - \frac{1}{2} \ep \sum_i  \frac{\partial H}{\partial \bq}(\bp_i,  \bq_i)  \frac{\partial}{\partial \bp_i}) (\bq_2, \bp_2)^t  .
\end{eqnarray}

And the substitution is given the algorithm for one time-step $n \rightarrow n+1$ and we obtain the solution in $t^{n+1}$:

$(\bq(t^{n+1}), \bv(t^{n+1}))^t = (\bq_3, \bv_3)^t$.

\end{algorithm12}

For studying the delicate higher accurate stability of the
Levitron, it is necessary to improve the standard St\"ormer-Verlet 
schemes, which are only second order scheme, but can be improve to a more accurate higher order scheme.
In the following, we discuss the extrapolation idea with respect to
the basic St\"ormer-Verlet algorithms.

\section{Improvement of the Splitting schemes with Extrapolation methods}
\label{extrapolation}

In the following, we discuss the different splitting schemes, that are based on the Strang-splitting scheme, see \cite{strang68} and extrapolated with the
so called MPE method, see \cite{chin_2011} for time-dependent problems.
In the first part we present the linear version, while in the second part we 
embed an iterative scheme to derive a nonlinear version.

The solution to the differential equation (\ref{splitt_1}) can
be formally written as  

\begin{equation}
u(t+\dt)=\T\Bigl(\exp\int_t^{t+\dt} {\mathcal A}(s)ds\Bigr) u(t),
\label{expth}
\end{equation}
where we assume general time-dependent operators ${\mathcal A}(t) = A(t) + B(t)$.

Following by Suzuki \cite{suzu93}, we have a 
{\it forward time derivative} operator, also called super-operator:
\begin{equation}
D={\frac{\buildrel \leftarrow\over\partial}{\partial t}}
\label{ftsh}
\end{equation}
such that for any two time-dependent functions $F(t)$ and $G(t)$,
\begin{equation}
F(t){\rm e}^{\dt D}G(t)=F(t+\dt)G(t).
\label{fg}
\end{equation}
If $F(t) = 1$, we have
\begin{equation}
1 {\rm e}^{\dt D}G(t)={\rm e}^{\dt D}G(t) = G(t).
\label{fg}
\end{equation}

By comparing with Trotters formula we can apply the Suzuki's decomposition of
the time-ordered exponential and obtain:
\begin{equation}
\T\Bigl(\exp\int_t^{t+\dt} {\mathcal A}(s)ds\Bigr)=\exp[\dt({\mathcal A}(t)+D)].
\label{tdecom}
\end{equation}
Thus time-ordering can be achieve by splitting an additional
operator $D$. 

With such a scheme, we can transforms in any existing splitting 
algorithms into integrators of non-autonomous equations. 

\begin{corrolar}
We achieve the following important second order 
symmetric splitting scheme for the St\"ormer-Verlet (PV) scheme, see also 
Figure \ref{verlet}:
\be
\T_2(\ep)=\e^{\frac12 \ep D}\e^{\ep {\mathcal A}(t)}\e^{\frac12 \ep D}=  \e^{\frac12 \ep B(t+\frac34 \ep)} \e^{\ep A(t+\frac12 \ep)} \e^{\frac12 \ep B(t+\frac14 \ep)} ,
\la{second}              
\ee
which is the second-order scheme with the assumption of
the commutation of the $A$ and $B$ between the $D$ operator, see also \cite{blanes2007}.
\end{corrolar}

\begin{proof}
For the second order algorithm, we apply the
Strang-splitting scheme for the three operators $A(t), B(t), D$
and we have assumed:  $[A(t), D] = 0$, $[B(t), D] = 0$.
\begin{eqnarray}
\T_2(\ep) && = \e^{\frac12 \ep D} \e^{\frac12 \ep B(t + \frac12 \ep D)}   \e^{\ep A(t)}  \e^{\frac12 \ep B(t)} \e^{\frac12 \ep D}  \nonumber \\
&& = \e^{\frac12 \ep B(t)}  \e^{\frac12 \ep D} \e^{\ep A(t+\frac12 \ep)}   \e^{\frac14 \ep D} \e^{\frac12 \ep B(t+\frac14 \ep)}  \e^{\frac14 \ep D}  \nonumber \\
&& = \e^{\frac12 \ep B(t)}  \e^{\frac14 \ep D} \e^{\ep A(t+\frac12 \ep)}   \e^{\frac12 \ep D} \e^{\frac12 \ep B(t+\frac14 \ep)}  \e^{\frac14 \ep D}  \nonumber \\
&& = \e^{\frac12 \ep B(t+\frac34 \ep)} \e^{\ep A(t+\frac12 \ep)} \e^{\frac12 \ep B(t+\frac14 \ep)} ,    
\end{eqnarray}
where we have applied the commutativity  with the $D$ operator and 
the shift with the forward time derivative operators.
\end{proof}

\begin{remark}
Every occurrence of the operator $\e^{d_i\ep D}$, from right to left, updates the current 
time $t$ to $t+d_i\ep$. If $t$ is the time at the start of the
algorithm, then after the first occurrence of $\e^{\frac12 \ep D}$, time is $t+\frac12\dt$.
After the second $\e^{\frac12 \ep D}$, time is $t+\dt$. Thus the leftmost 
$\e^{\frac12 \ep D}$ is not without effect, it correctly updates the time for the next
iteration, see also \cite{suzu93}. 

Thus the iterations of $\T_2(\ep)$ implicitly imply
\ba
 \T_2^2(\ep/2)&=&\e^{\frac12\ep {\mathcal A}(t+\frac34 \ep)}\e^{\frac12\ep {\mathcal A}(t+\frac14 \ep)} ,
\la{tn}
\ea
by inserting equation (\ref{second}), we obtained:
\begin{eqnarray}
&& \T_2^2(\ep/2)\\
&=&\e^{\frac14 \ep B(t+\frac78 \ep)} \e^{\frac12 \ep A(t+\frac34 \ep)} \e^{\frac14 \ep B(t+\frac58 \ep)} \e^{\frac14 \ep B(t+\frac38 \ep)} \e^{\frac12 \ep A(t+\frac14 \ep)} \e^{\frac14 \ep B(t+\frac18 \ep)} . \nn
\end{eqnarray}

\end{remark}

For higher orders we have explicitly:
\be
\T_4(\ep)=-\frac13\T_2(\ep)
+\frac43\T_2^2\left(\frac\ep{2}\right)
\la{four} ,
\ee
\be
\T_6(\ep)=\frac1{24} \T_2(\ep)
-\frac{16}{15}\T_2^2\left(\frac\ep{2}\right)
+\frac{81}{40}\T_2^3\left(\frac\ep{3}\right)
\la{six} ,
\ee
\be
\T_8(\ep)=-\frac1{360} \T_2(\ep)
+\frac{16}{45}\T_2^2\left(\frac\ep{2}\right)
-\frac{729}{280}\T_2^3\left(\frac\ep{3}\right)
+\frac{1024}{315}\T_2^4\left(\frac\ep{4}\right)
\la{eight} ,
\ee
\ba
&&\T_{10}(\ep)=\frac1{8640} \T_2(\ep)
-\frac{64}{945}\T_2^2\left(\frac\ep{2}\right)
+\frac{6561}{4480}\T_2^3\left(\frac\ep{3}\right)\nn\\
&&\qquad\qquad\quad-\frac{16384}{2835}\T_2^4\left(\frac\ep{4}\right)
+\frac{390625}{72576}\T_2^5\left(\frac\ep{5}\right).
\la{ten}
\ea

\begin{remark} In the work of  Blanes, Casas and Ros\cite{bcr99} and Chan and Murua\cite{cm00} the idea of extrapolating symplectic algorithms are also discussed. They presented the case of extrapolating an $2n$-order symplectic integrators and noted that extrapolating a $2n$-order symplectic integrator will preserve the symplectic character 
of the algorithm to order $4n+1$. 
\end{remark}

\section{Numerical Results}
\label{numerical}

In the following we deal with the computation of the 
Hamiltonian, which is dervied in Section \ref{model},
see also  \cite{dull98,gans97}.

Our Hamiltonian of the Levitron is given as:
\begin{eqnarray}
\label{equ1_ham}
\nonumber H=&\frac{1}{2}\left(p_1^2+p_2^2+p_3^2+\frac{p_4^2}{a}+\frac{(p_5-p_6 \cos q_4)^2}{a \sin ^2 q_4} + \frac{p_6}{c} \right) \\
&- M \left[\sin q_4 \left( \cos q_5 \frac{\pa \Psi}{\pa q_1} \sin q_5 \frac{\pa \Psi}{\pa q_2} \right) +\cos q_4 \frac{\pa \Psi}{\pa q_3}  \right] +q_3
\end{eqnarray}

For our splitting scheme, we apply the Hamiltonian of (\ref{equ1_ham}), we have:
\begin{eqnarray}
\label{decouple_1}
{\bf \dot q} & =& \frac{\partial H}{\partial \bp}(\bp,  \bq) \\
& = & \left( p_1, p_2, p_3, \frac{p_4}{a}, \frac{\left(p_5 - p_6 \cos q_4 \right)^2}{a \sin^2 q_4}, \frac{p_6 (\cos^2 q_4 + (a /c) \sin^2 _4) - p_5 \cos q4}{a \sin^2 q_4} \right) \nonumber 
\label{ham}
\end{eqnarray}
given as operator $A$ and
\begin{align}
\nonumber {\bf \dot p} =& - \frac{\partial H}{\partial \bq}(\bp, \bq) \\
\nonumber &  = ( M\left( \sin q_4 \cos q_5 \frac{\pa^2 \Psi}{\pa q_1^2} + \cos q_4 \frac{\pa^2 \Psi}{\pa q_1 \pa q_3}  \right), M\left( \sin q_4 \cos q_5 \frac{\pa^2 \Psi}{\pa q_2^2} + \cos q_4 \frac{\pa^2 \Psi}{\pa q_2 \pa q_3}  \right), \\
\nonumber &  M\left( \sin q_4 \left(\sin q_5 \frac{\pa^2 \Psi}{\pa q_2 \pa q_3} + \cos q_5 \frac{\pa^2 \Psi}{\pa q_1 \pa q_3} \right) + \cos q_4 \frac{\pa^2 \Psi}{\pa q_3^2} \right)-1, \\
\nonumber &  M\left( \cos q_4 \left( \sin q_5 \frac{\Psi}{q_2} + \cos q_5 \frac{\pa \Psi}{\pa q_1} \right) - \sin q_4 \frac{\pa \Psi}{\pa q_3} \right) \\
& - \frac{p_6(p_5-p_6 \cos q_4)}{a \sin q_4}-\frac{\cos q_4(p_5-p_6\cos q_4)^2}{a \sin ^3 q_4}, \nonumber \\
 &  M\left( \sin q_4 \left( \cos q_5 \frac{\pa \Psi}{\pa q_2} - \sin q_5 \frac{\pa \Psi}{\pa q_1} \right) \right),
0  )
\label{decouple_2}
\end{align}
given as operator $B$, which we insert into the Algorithm \ref{algo_2}.

\begin{figure}[htbp]
 \begin{center}	
	\includegraphics[width=0.3\textwidth,angle=-90]{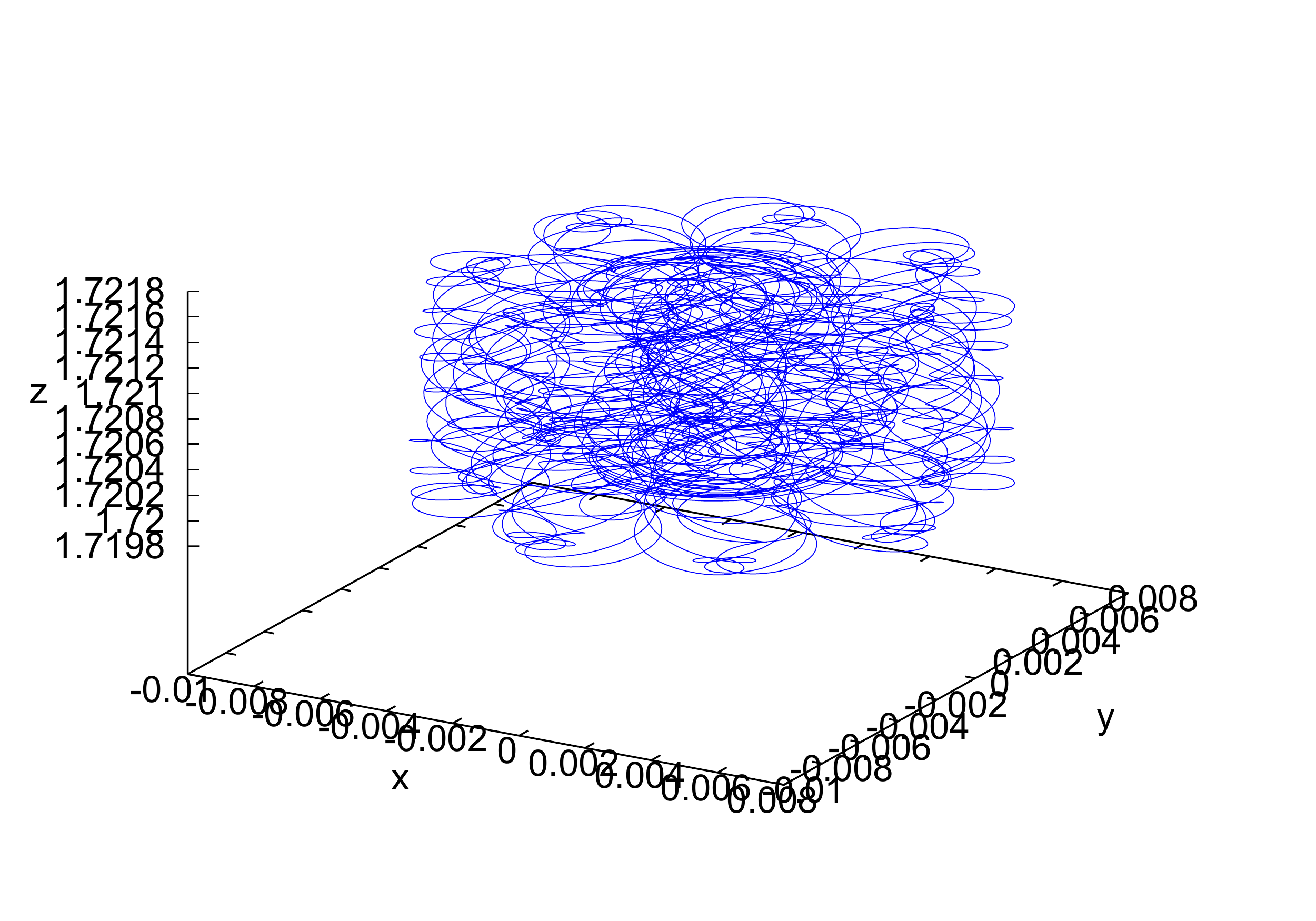}
%	\centering
%	\includegraphics[width=\textwidth]{images/xy_verlet.eps}
	\includegraphics[width=0.3\textwidth,angle=-90]{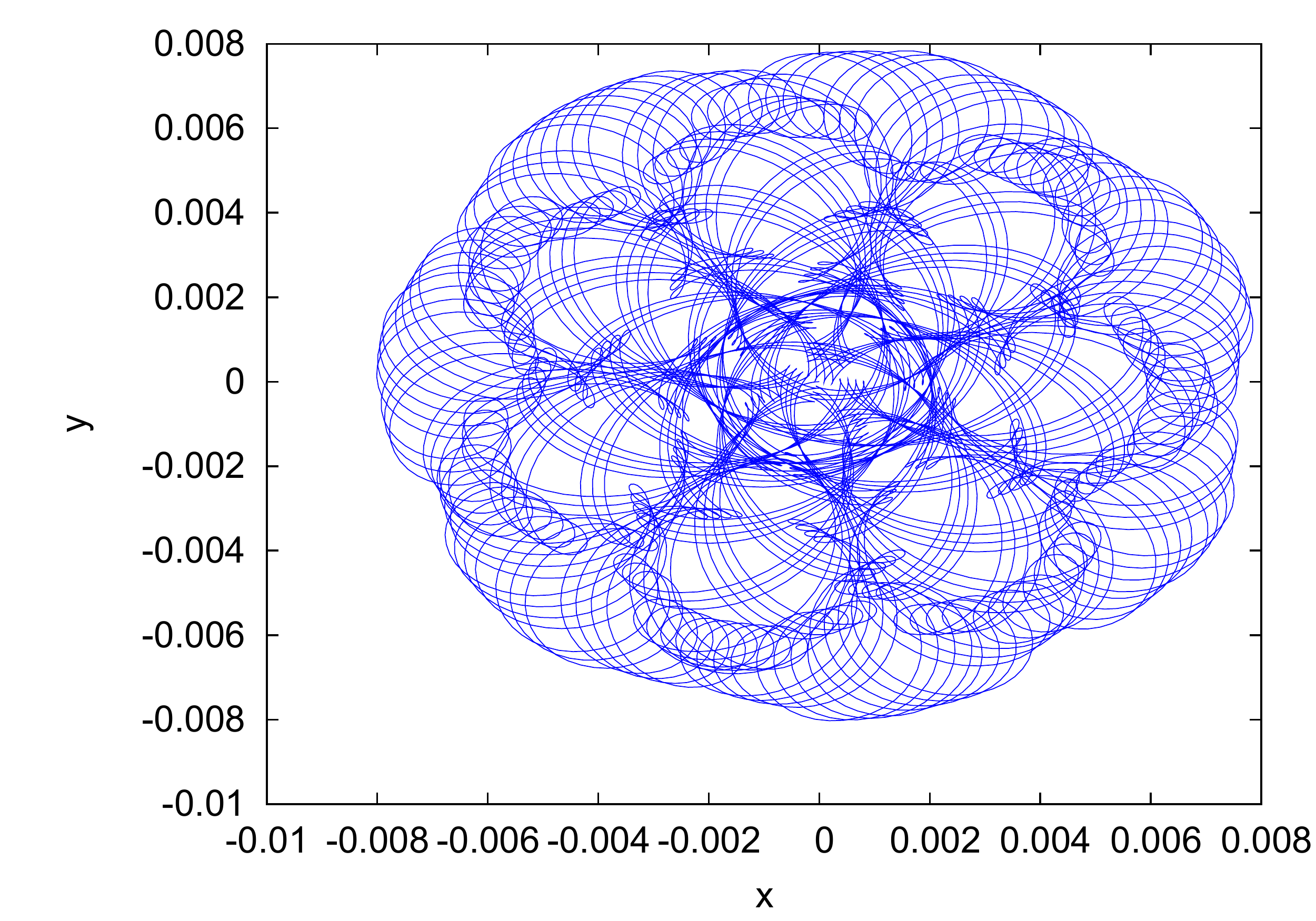}
 \caption{Trajectory calculated with Verlet algorithm (Left figure: 3D presentation, right figure: 2D presentation).} 
 \label{fig:trajectory_verlet}
\end{center}
\end{figure}

We compare our novel schemes (extrapolated St\"omer-Verlet method) 
with standard and Runge-Kutta algorithms.
Due to the long computation time needed, we simulated only 1000 timesteps and compare the trajectory with the reference solution from the Runge-Kutta algorithm. In figure \ref{fig:trajectory_verlet} is shown how the trajectory of the same initial conditions looks like with the Verlet algorithm. 

\begin{figure}[htbp]
\centering	
	\begin{minipage}[b]{0.47\textwidth}
	\centering 
	\includegraphics[width=\textwidth]{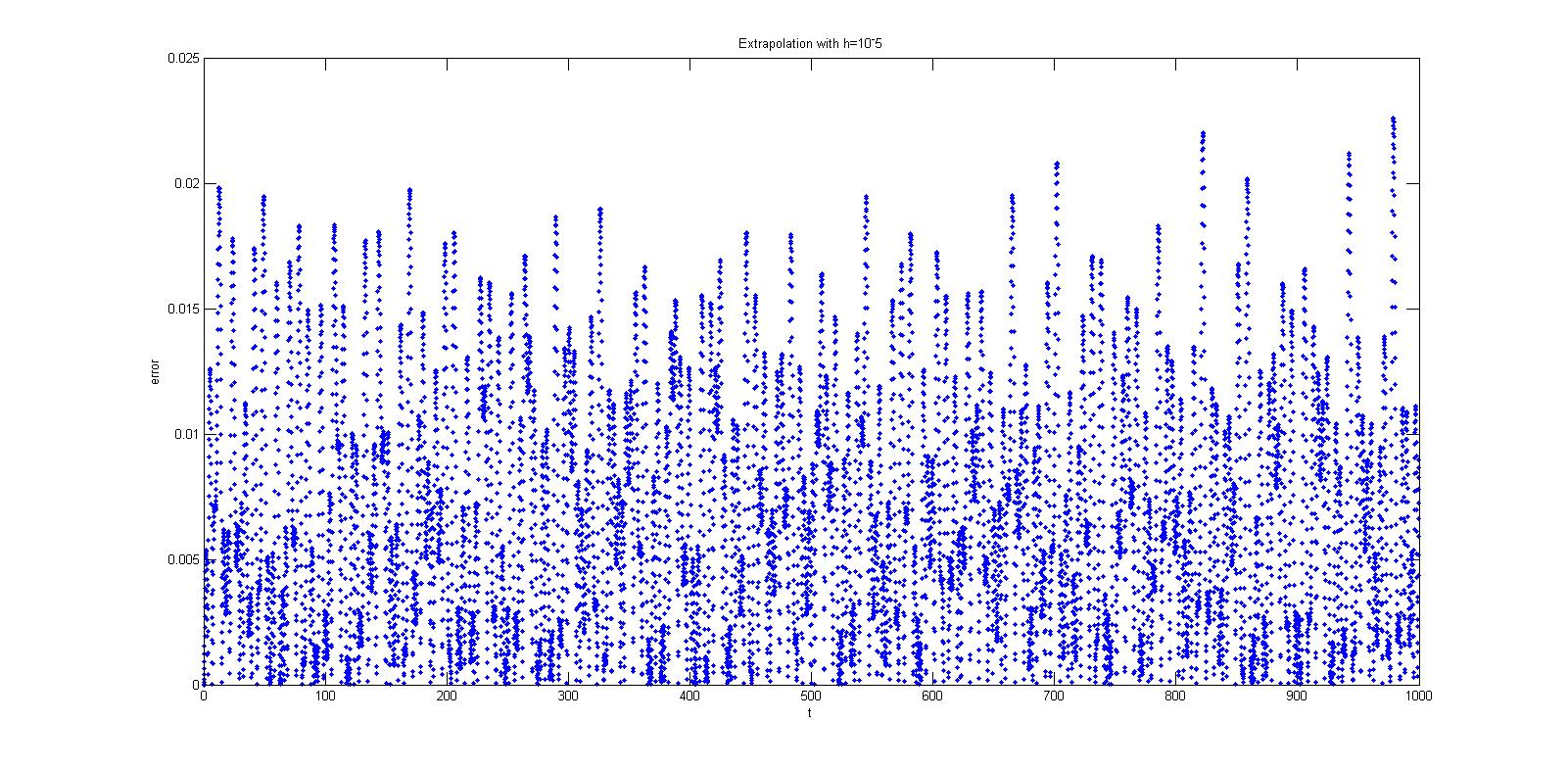}
	\end{minipage}
	\hfill
	\begin{minipage}[b]{0.47\textwidth}
	\centering 
	\includegraphics[width=\textwidth]{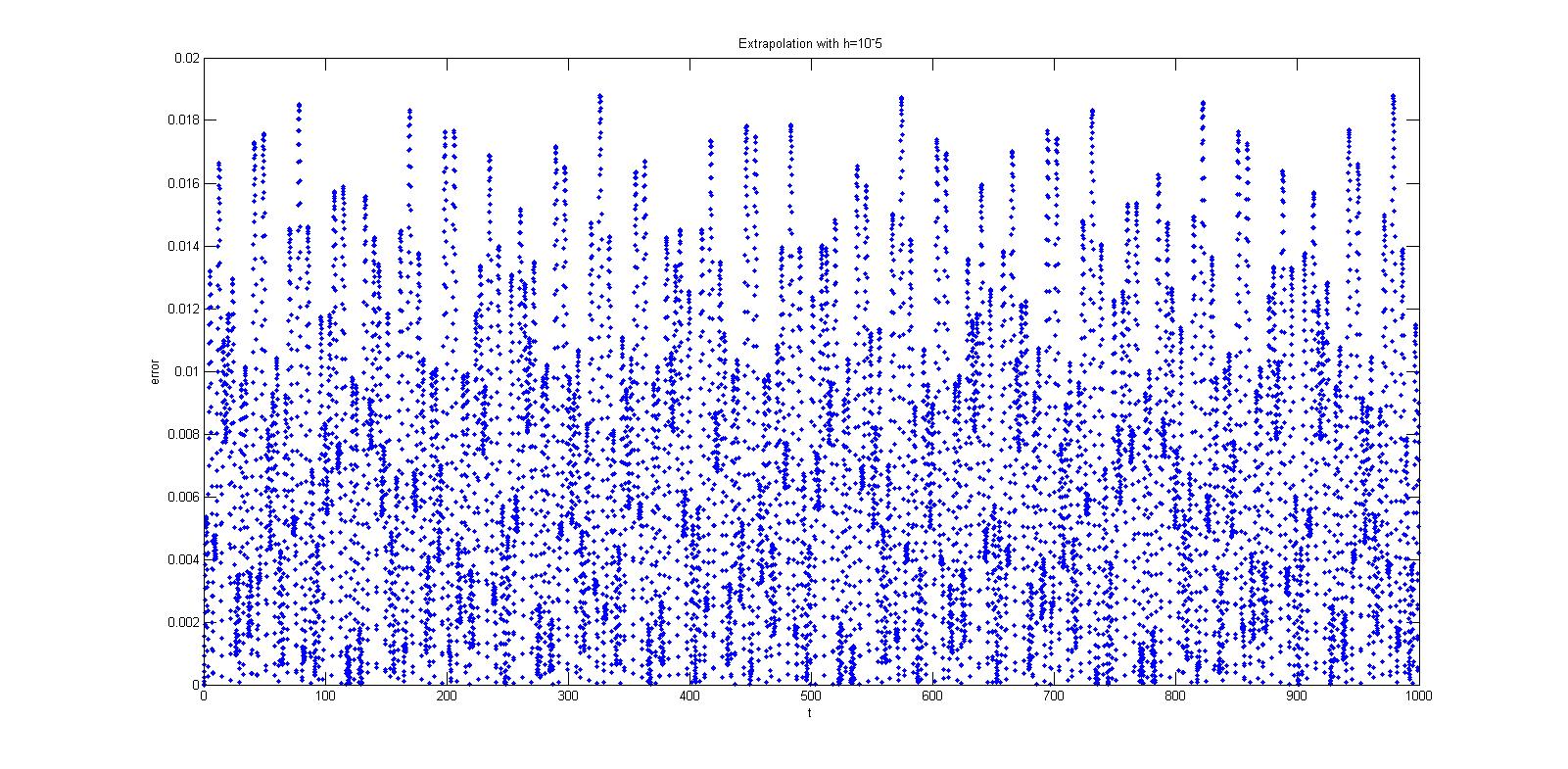}
	\end{minipage}
	\caption{Errors of the numerical scheme: Extrapolation Scheme with Verlet method a Kernel (left figure: 4th order scheme with $h=10^{-5}$ and 6th order scheme with $h=10^{-6}$).} 
 \label{fig:extra_1}
\end{figure}

\begin{figure}[htbp]
\centering	
	\begin{minipage}[b]{0.47\textwidth}
	\centering 
	\includegraphics[width=\textwidth]{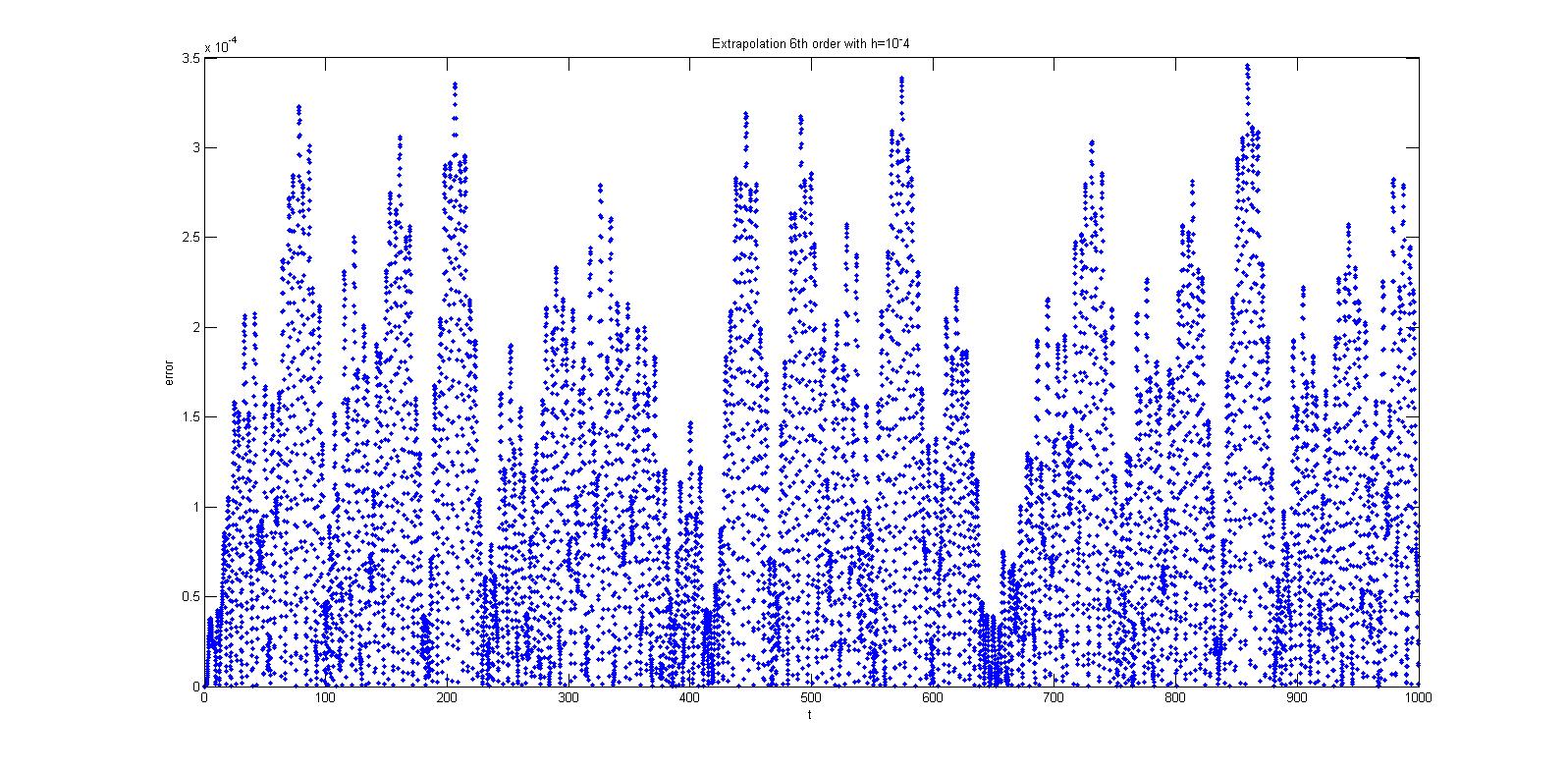}
	\end{minipage}
	\hfill
	\begin{minipage}[b]{0.47\textwidth}
	\centering
	\includegraphics[width=\textwidth]{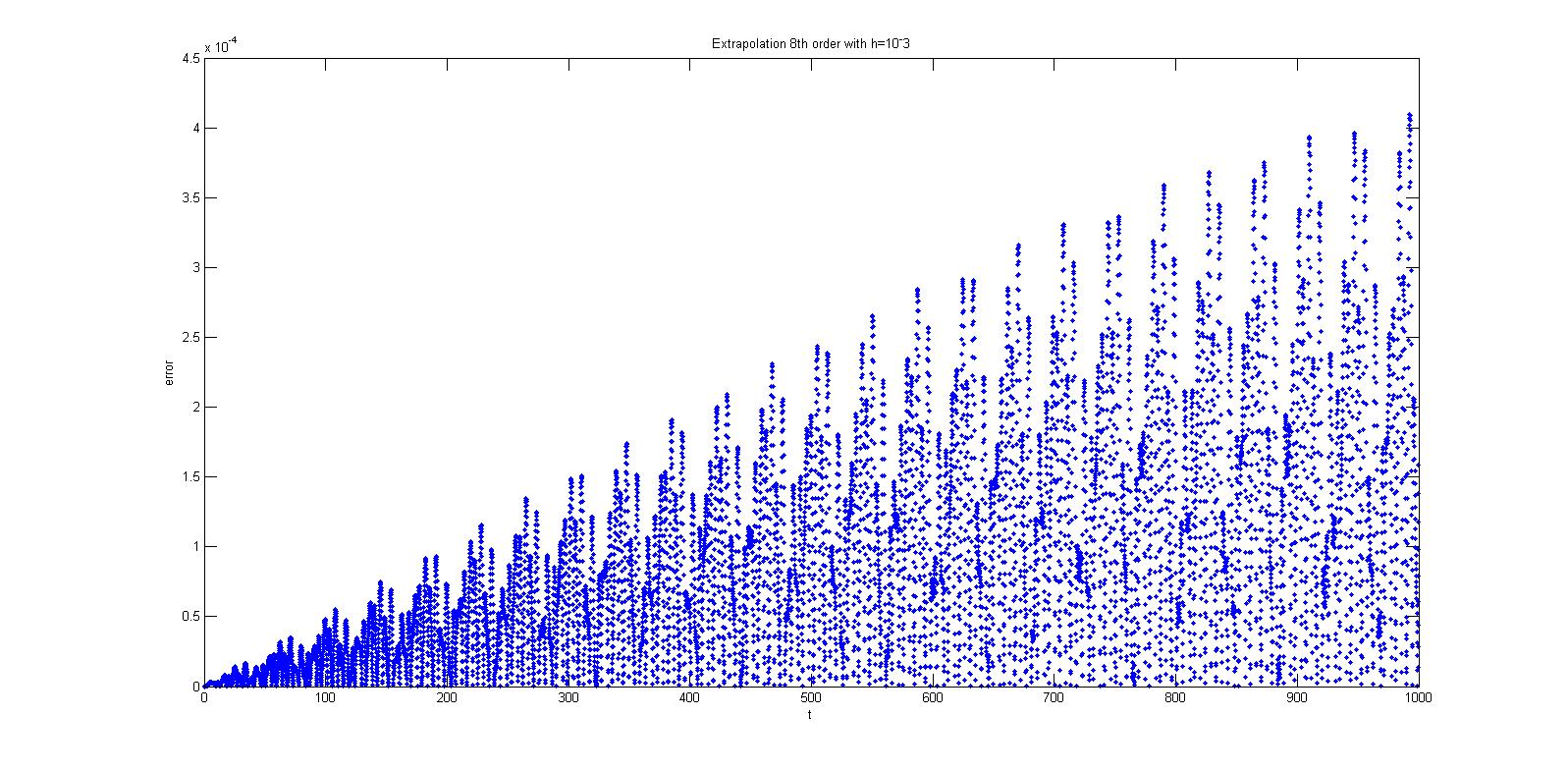}
	\end{minipage}
	\caption{Errors of the numerical scheme: Extrapolation Scheme with Verlet method a Kernel (left figure: 6th order scheme with $h=10^{-4}$ and 8th order scheme with $h=10^{-3}$).} 
 \label{fig:extra_2}
\end{figure}

We improve the solution with an extrapolation scheme with fourth order. 
We have a view at the errors this algorithm produces in comparison with the Runge-Kutta Solution with small time-steps ($10^{-5}$ time units per step).
In Figure \ref{fig:extra_1} and \ref{fig:extra_2}, we presented the results of the 4th, 6th and 8th order
Multiproduct expansion method with different time-steps and compared it with a fine resolved 4th order Runge-Kutta Benchmark solution ($h = 10^{-8}$).

The time scales and computational amount for the extrapolation schemes are given in Table \ref{table3} and \ref{table4}.
\begin{table}[ht]
\centering
\begin{tabular}{|c|c|c|c|c|}
\hline  & \multicolumn{2}{|c|}{Extrapolation 4th order} & \multicolumn{2}{|c|}{Extrapolation 6th order}  \\
\hline timestep &  $10^{-5}$ & $10^{-6}$ & $10^{-5}$ & $10^{-6}$   \\ 
\hline number of steps & $10^{8}$ & $10^{9}$ & $10^{8}$ & $10^{9}$ \\ 
\hline computing time & 14min & 142min & 29min & 272min  \\
\hline mean error & 0.007  & 0.007 & 0.0068 & 0.0068\\
\hline maximal error & 0.0226 & 0.0234 & 0.0188 & 0.0188\\
\hline
\end{tabular} 
\caption{Errors and Computational Time with 4th order MPE scheme using
Verlet Scheme as Kernel.}
\label{table3}
\end{table}

\begin{table}[ht]
\centering
\begin{tabular}{|c|c|c|c|c|}
\hline 
 & Extrapolation & Extrapolation & \multicolumn{2}{|c|}{Extrapolation}  \\
 & 6th order & 8th order & \multicolumn{2}{|c|}{10th order}  \\
\hline timestep &  $10^{-4}$ & $10^{-3}$ & $10^{-2}$ & $10^{-3}$  \\ 
\hline number of steps & $10^{7}$ & $10^{6}$ & $10^{5}$  & $10^{6}$ \\ 
\hline computing time & 2.5min & 0.5min & 3sec & 32.8sec \\
\hline mean error & $1.0244\cdot 10^{-4}$  & $9.6297\cdot 10^{-5}$ & $0.013$ & $2.2397\cdot 10^{-5}$ \\
\hline maximal error & $3.4608\cdot 10^{-4}$ & $4.0936\cdot 10^{-4}$ & $0.01$ & $9.8868\cdot 10^{-5}$ \\
\hline
\end{tabular} 
\caption{Errors and Computational Time with higher order MPE scheme using
Verlet Scheme as Kernel.}
\label{table4}
\end{table}

We have also tested the 10th order extrapolation with $10^{-2}$ time units 
per step and also obtained stable trajectory.

\begin{remark}
In the examples, we have verified, that we can improve a
 basic second order symplectic splitting scheme with extrapolation schemes. 
At least achieved higher accurate solutions and save computational time.
Moreover we save computer resources and obtained stable trajectories with 
larger time-steps.
The best result we achieve with the
order 10 and $h=10^{-2}$ for such a case we could improve the results and
are $10$-times faster than with standard 4th order explicit Runge-Kutta schemes.
\end{remark}

\section{Conclusion}
\label{conclusion}

In the paper, we have presented a model to simulate a Levitron.
Based on the given Hamiltonian system, which can be written as 
large system of time-dependent ordinary differential equation, we 
present novel and faster solvers based on splitting and extrapolation ideas. 
We could achieve more accurate and stable
results with higher order schemes and save computational time with 
respect of stable computations.
In future, we concentrate on the numerical analysis and 
embedding higher order splitting kernels to nonlinear differential
equations based on Hamiltonian systems.


\begin{thebibliography}{10}

\bibliographystyle{plain}

% \bibitem[Names(Year)]{label} or \bibitem[Names(Year)Long names]{label}.
% (\harvarditem{Name}{Year}{label} is also supported.)

\bibitem{geiser_2012} 
Geiser J,  L\"uskow K. Splitting methods for Levitron Problems, Preprint, Department of Physics, EMU University of Greifswald, Germany, 2011.



\bibitem{dullin2004}
Dullin H R. Poisson Integrator For Symmetric Rigid Bodies. Regular and chaotic dynamics. 9(3) (2004) 255-264.


\bibitem{chin_2011} 
Chin S, Geiser J. Multi-product operator splitting as a general method of solving autonomous and non-autonomous equations, IMA Journal of Numerical Analysis, 31 (2011) 1552-1577.



\bibitem{dull98}
Dullin H R, Easton R. Stability of Levitron. Physica D: Nonlinear Phenomena, 126(1-2) (1999) 1-17.



\bibitem{gans97}
Gans R F, Jones T B, Washizu M. Dynamics of the Levitron. J. Phys. D. 31(1998) 671-679, 1998.

\bibitem{hair2003}
Hairer E, Lubich Chr,Wanner G. Geometric numerical integration illustrated by the St\"ormer-Verlet method. Acta Numerica (2003) 399-450.



\bibitem{goldstein81}
Goldstein H, Poole Ch P, Safko J. Classical mechanics. Addison Wesley, San Francisco, USA, 2002.





 \bibitem{strang68} 
Strang G. On the construction and comparison of difference schemes. SIAM J. Numer. Anal. 5 (1968) 506-517.



\bibitem{suzu93}
Suzuki M. General Decomposition Theory of Ordered Exponentials.
Proc. Japan Acad., 69B (1993) 161-166.

\bibitem{blanes2007} 
Blanes S, Casas F, Murua, A.  Splitting method for non-autonomous 
linear systems. J. Comput. Math. 84 (2007) 713-727.



\bibitem{bcr99} 
Blanes S, Casas F, Ros J.  Extrapolation of symplectic integrators.  Celest. Mech. Dyn. Astron.  75 (1999) 149-161.


\bibitem{cm00} Chan R, Murus A. Extrapolation of symplectic methods for Hamiltonian problems. Appl. Numer. Math. 34(2000) 189-205.

\end{thebibliography}
\end{document}